\documentstyle[12pt,fleqn,twoside]{article}
\def\R{\textup{I$\!$R}}

\newtheorem{thm}{Theorem}[section]  %% Definition of Theorem
\newtheorem{crlr}[thm]{Corollary}      %% Definition of Corollary
\newtheorem{prp}{Proposition }[section]     %% Definition of Proposition
\newtheorem{prp4.}{Proposition 4.}%[section]     %% Definition of Proposition%\begin{prp}
\newtheorem{Def}{Definition }[section]    %% Definition of Definition
\newtheorem{rem}[thm]{Remark}         %% Remark
\newtheorem{dim.}{Proof.}     %% Definition of Proposition
\makeatletter
\def\theequation{\thesection.\@arabic\c@equation}
\def\thethm{\thesection.\@arabic\c@thm}
\def\thelem{\thesection.\@arabic\c@thm}
\def\thecrlr{\thesection.\@arabic\c@thm}
\def\theprp{\thesection.\@arabic\c@thm}
\def\therem{\thesection.\@arabic\c@thm}
\makeatother

\date{}

\begin{document}

\baselineskip=9mm
\title{Mixed Morrey spaces and their applications\\
to partial differential equations }
%\author{{\sc Maria Alessandra Ragusa and Andrea Scapellato}
%\thanks{ Dipartimento di Matematica e Informatica, Universit\`{a} di Catania,
%Viale Andrea Doria 6, 95125 Catania (Italy).
%Ccorresponding author M.A.Ragusa, e-mail: maragusa@dmi.unict.it
%% \noindent { Acknowledgments.}
%%Acknowledgments. The  author is greatly indebted to the anonymous referee
%%for the helpful ideas and suggestions.
%$$\,$$
%Mathematics Subject Classification (2000): Primary 31B10, 43A15, 35K20. Secondary
%32A37, 46E35.\qquad \qquad \qquad \, Key words and phrases: Morrey spaces, mixed norm, regularity results, partial differential equations.
%}}
\author{
Maria Alessandra Ragusa
\thanks{Dipartimento di Matematica e Informatica, Universit\`{a} di Catania,
Viale Andrea Doria, 6-95125 Catania, Italy,
e-mail:maragusa@dmi.unict.it, corresponding author.}
~\&~ Andrea Scapellato
\thanks{Dipartimento di Matematica e Informatica, Universit\`{a} di Catania,
Viale Andrea Doria, 6-95125 Catania, Italy,
e-mail:scapellato@dmi.unict.it}
%\thanks{%This work was partly supported by FIRB 2014
%Mathematics Subject Classification (2000): Primary 31B10, 43A15, 35K20. Secondary
%32A37, 46E35.\qquad \qquad \qquad \, Key words and phrases: Morrey spaces, mixed norm, regularity results, partial differential equations
%}
}
\date{ }

\maketitle

\setcounter{page}{1}
\setcounter{equation}{0}
\setcounter{thm}{0}

\centerline{ Abstract}

In this paper, new classes of functions are defined.
These spaces generalize Morrey spaces and give a refinement of
Lebesgue spaces.  Some  embeddings  between these new classes are also proved.
Finally, the authors apply these classes of functions to obtain regularity results for solutions of partial differential equations of %elliptic
parabolic type.
%------------------------------------------------------------------

\section{Introduction} %Section 1
\setcounter{equation}{0}
\setcounter{thm}{0}

%The purpose of this note is to define %${\cal R}^{(p,q,\lambda)}$
This paper aims at defining new spaces and to study some embeddings  between them.  We will refer to them with the symbol $ L^{q,\mu}(0,T,L^{p,\lambda}(%\Omega
\R^n)). $ As applications we obtain some estimates, in these classes of functions, for the solutions of partial differential equations of %elliptic
 parabolic type in nondivergence form. Preparatory to achieving these results is the study %, in these classes,
 of the behaviour of Hardy-Littlewood Maximal function, Riesz potential, Sharp and Fractional maximal functions, Singular integral operators with Calder\'on-Zygmund kernel and Commutators (see e.g. %\cite{DR2},
 \cite{RDUKE} \cite{JOTA}).

We stress that are obtained results, known in $L^p, $ in a new class of functions that can be view as an extension of the Morrey class introduced in 1966 in \cite{Morrey}, and used by a lot of authors, see e.g. in \cite{Caffarelli}, %\cite{PS},
 recently in \cite{RTT}, \cite{PRSW}, \cite{G1}, \cite{G2}, \cite{G3} and others.

Let us point out that in doing this we need an extension to $ L^{q,\mu}(0,T,L^{p,\lambda}(%\Omega
\R^n)) $ of a celebrated inequality of Fefferman and Stein (see \cite{FeStein}) concerning the Sharp and the Maximal function (%Proposition
Theorem \ref{%prp3
thmprp1}) and, also, we study the behaviour of Riesz potential in the new class of functions, obtaining an extension of both a known estimate originally proved by Adams in \cite{Adams} %and %another result %by Spanne
as well as of a result announced by Peetre in \cite{Peetre}.

\section{Definitions and Preliminary Tools} %Section 2
\setcounter{equation}{0}
\setcounter{thm}{0}

In the sequel let $T>0 $ and $\Omega $ be a bounded open set of $\R^{n} $ such that $\exists A>0\, : \forall x \in \Omega$ and $ 0 \leq \rho \leq \,diam \,(\Omega), \,$  $|Q(x,\rho) \cap \Omega| \geq A\, \rho^n, $ being $Q(x,\rho) $ a cube centered in $x, $ having edges parallel to the coordinate axes and lenght $2 \rho$.

\begin{Def}
Let $1\!<\!p\!<\!+ \infty, $ $0\!<\!\lambda\!<\!n $ and $f $ be a real measurable
function defined in %on the open bounded set
 $\Omega\,\subset\, \R^n.$
%(n \geq 1).$

If $|f|^p\,$ is summable in $\Omega$  and the set described by the quantity
\begin{equation}\label{1}
\frac{1}{\rho^\lambda} \int_{\Omega \cap B_\rho (x)} |f(y)|^p\,dy,
\end{equation}
when changing of  $\rho $  in $]0, \,diam\,\Omega[ $ and $x \in \Omega,
$ has an upper bound, then we say that $f $ belongs to the {\it Morrey Space }
$L^{p,\lambda}(\Omega). $
\end{Def}

If $f \in L^{p,\lambda}(\Omega), $ we define
\begin{equation}\label{2}
\|f\|^p_{L^{p,\lambda}(\Omega)} := \sup_{%\substack
{x \in \Omega \atop  % 3 giugno 2013
%\
\ %0<\rho< diam\, \Omega
\rho> 0}} \,\frac{1}{\rho^\lambda}
\int_{\Omega \cap B_\rho (x)} |f(y)|^p\,dy
\end{equation}
and the vector space naturally associated to the set of functions
 in $L^p(\Omega) $ such that (\ref{2}) is finite, endowed with the
norm (\ref{2}), is a normed and complete space.\\
The exponent $\lambda $  can take values that are not
belonging to $]0, n[ $ but the unique cases of real interest are
that one for which $\lambda \in ]0, n[. $

The above defined spaces are used, among others, in the theory of regular
solutions to nonlinear %elliptic
 partial differential equations and for the study of
local behavior of solutions to nonlinear %elliptic
 equations and
systems (see  e.g. \cite{Morrey}, \cite{MUSTA}).

\begin{rem}\label{rem2.1.}\hspace*{-0.6em}\textbf{.}
Similarly we can define the Morrey space in $L^{p,\lambda}(\R^n) $ as the space of functions such that is finite:
\begin{equation}\label{3}
\|f\|^p_{L^{p,\lambda}(\R^n)} := \sup_{%\substack
{x \in \R^n \atop
 \rho\,>\,0}} \,\frac{1}{\rho^\lambda}
\int_{ B_\rho (x)} |f(y)|^p\,dy.
\end{equation}
\end{rem}

\begin{Def} \label{def2.1.} \hspace*{-0.6em}\textbf{.} %${\cal L}(0,T,{\cal L}^{p,\lambda}(\Omega)). $
Let $1<p, q <+ \infty, $ $0<\lambda, \mu <n. $ We define the set
$L^{q,\mu}(0,T, L^{p,\lambda}(\Omega)) $ as the class of functions $f$ such that is finite:
\begin{equation}\label{4}
\!\!\!\!\!\!\!\!
\!\!\!\!\!\!\!\!\|f\|_{L^{q,\mu}(0,T,L^{p,\lambda}(\Omega))} :=
\left(
\sup_{%\substack
{t_0, t \in (0,T)\,\atop  \rho>0 }
}\,
\frac{1}{\rho^\mu}
\int_{(0,T) \cap (t_0-\rho, t_0+\rho)}
\left(
\sup_{%\substack
{x \in \Omega
\atop  %0<\rho< diam\, \Omega
\rho > 0}} \,\frac{1}{\rho^\lambda}
\int_{\Omega \cap B_\rho (x)} |f(y,t)|^p\,dy
\right)^{\!\!\frac{q}{p}}
dt\,
\right)^{\!\!\frac{1}{q}}\!\!\!,
\end{equation}
with obvious modifications if $\Omega = \R^n. $
\end{Def}

\begin{Def}\label{def2.6.}
Let $f$ be a locally integrable function defined on ${\R}^{n}.$ We say
 that $f$ is in the space $BMO( {\R}^{n} )$ (see \cite{JN}) if
$$
%\| f \|_* \equiv
\sup_{B \subset {\R}^{n}} \frac{1}{|B|}
\int\limits_B |f(y)-f_B| dy < \infty
$$
where $B$ runs over the class of all balls %parabolic cubes
 in ${\R}^{n} $ and
$f_B%Q
 = \frac{1}{|B%Q
 |} \int\limits_B%Q
  f(y)dy.$
\end{Def}

Let $f\in BMO( {\R}^{n})$  and $r>0.$ We define the $VMO$ modulus
of $f$ by the rule
$$
\eta (r) = \sup_{
%{x \in {\R}^n,}\
{\rho \leq r}}
\frac{1}{|B%Q
_\rho|} \int\limits_{B%Q
_\rho} |f(y)-f_{B%Q
_\rho}| dy
$$
where $B%Q
_{\rho}$ is a %parabolic cube
 ball  with radius $\rho ,$ $\rho \leq r.$

$BMO$ is a Banach space with the norm $\| f \|_*\,=\,\sup_{r > 0} \eta(r).$

\begin{Def}\label{def2.7.}
We say that a function $f\! \in \! \! BMO$ is in the Sarason class
 $VMO({\! \R}^{n}\!)$ (see \cite{S}) if
$$
\lim_{r \to 0^+} \eta (r) =0.
$$
\end{Def}

%Dalla section 5 fino a...
%\begin{Def}\label{Defmetric}
%Let  $x,y \in \R^{n} $  %$x=(x',t)=(x'_1, x^_2, \ldots, x'n, t), $
% and endow $\R^{n} $ with the known parabolic metric, introduced by Fabes and Rivi\'ere in the paper \cite{FR},
%$$%\begin{equation}
%d(x,y)=\rho (x-y),
%$$%\end{equation}
%being
%\begin{equation}
% \rho(x)(x)\,=\, \sqrt{ \frac{|x'|^2+\sqrt{|x'|^4 + 4 t^2} }{2}}, \qquad x=(x',t)\,=\,(x'_1, \ldots,x'_{n-1},t)\in \R^{n}.\,
%\end{equation}
%
%\end{Def}
\begin{Def}\label{Def5.1}%[pala Soft Pot An 2004, pg.240]
Let %$\alpha_1,\ldots,\alpha_n $ be real numbers, $\alpha_i\geq 0,\,$ $\alpha\,=\,\sum_{i=1}^n \alpha_i, \,$
$\Sigma %S^{n-1}
\,$ the unit sphere:
$\Sigma%S^{n-1}\,
=\,\{x \in \R^{n+1},\,|x|\,=\,1\}.$

The function $k: \R^{n+1} \backslash \{0\} \to \R\,$ %is  {\it kernel %with mixed homogeneity}
is the classical Calder\'on-Zygmund kernel  if:

$1)\, k \in C^\infty(\R^{n+1} \backslash \{0\});$

$2) k(\mu%^{\alpha_1}
x_1, \mu%^{\alpha_2}
x_2, \ldots,\mu%^{\alpha_n}
x_n, \mu^2 t)\,=\,\mu^{-(n+2)%\alpha
}
\,k(x),\,$ for each $ \mu>0;$

$3) \int_{%S^{n-1}
\Sigma}|k(x)|\,d\,\sigma_x < \infty\,\,$ and $\qquad\int_{%S^{n-1}
\Sigma}k(x) \,d \,\sigma_x =0. $
\end{Def}
%is the classical Calder\'on-Zygmund kernel.

%DA QUI C'E' IL NUCLEO VARIABILE COME IN PALA-SOFT POT.An 2004 PG:240

%\begin{Def}\label{Def5.1}%[pala Soft Pot An 2004, pg.240]
%Let $\alpha_1,\ldots,\alpha_n $ be real numbers, $\alpha_i\geq 0,\,$ $\alpha\,=\,\sum_{i=1}^n \alpha_i, \,$ $S^{n-1}\,$ is the unit shpere:
%$$S^{n-1}\,=\,\{x \in \R^n,\,|x|\,=\,1\}.$$
%
%A function $k: \R^n \backslash \{0\} \to \R\,$ is a {\it kernel with mixed homogeneity} if:
%
%$1)\, k \in C^\infty(\R^n \backslash \{0\});$
%
%$2) k(\mu^{\alpha_1}x_1, \mu^{\alpha_2}x_2, \ldots,\mu^{\alpha_n}x_n)\,=\,\mu^{-\alpha}\,k(x),\,$ for each $ \mu>0;$
%
%$3) \int_{S^{n-1}}|k(x)|\,d\,\sigma_x < \infty\,\,$ and $\qquad\int_{S^{n-1}}k(x) \,d \,\sigma_x =0. $
%\end{Def}
%
%Let us note that the special case $\alpha_i=1\,$ for $i=1,\ldots,n $ and, thus $\alpha=n,\,$ gives the classical Calder\'on-Zygmund kernel.

The above definition, in particular condition $2), $ suggest to endow $\R^{n+1} $ with a metric, different to the standard Euclidean one. Thus let us consider, as Fabes and Rivi\'ere in the celebrated paper \cite{FR}, the following distance $d(x,y)=\rho(x-y) $ between two generic points $x, y \in \R^{+1}n $ (used e.g. in \cite{BC}%pg.1736
),
\begin{equation}
\rho(x)\,=\, \sqrt{ \frac{|x'|^2+\sqrt{|x'|^4 + 4 t^2} }{2}}, \qquad x=(x',t)\,=\,(x'_1, \ldots,x'_{n},t)\in \R^{n+1},\,
\end{equation}
 Then $\R^{n+1},$ endowed with this metric, is a metric space.

%FINO A QUI C'E' IL NUCLEO VARIABILE COME IN PALA-SOFT POT.An 2004 PG:240

\begin{Def}\label{Def5.2}%[Pala Soft Pot An 2004, pg.241] [BC] pg.1747
 The function $k(x,y): \R^{n+1}\times \R^{n+1}\backslash\{0\}\to \R\,$ is a {\it variable } Calder\'on-Zygmund kernel %with mixed homogeneity
 if:

1) $k(x,\cdot)$ is a kernel in the sense of the above Definition \ref{Def5.1}, for a.e. $x \in \R^{n+1}$

2) $\sup_{\rho(y)=1} \left|  \left(\frac{\partial}{\partial y}\right)^\beta k(x,y) \right| \leq c(\beta), $ for every multi-index $\beta$, independently of $x.$

%where $\rho$ is such that $\R^n, $ endowed with the metric $\rho, $ is a metric space.

\end{Def}

%...fino a qua la sectio 5

Next Proposition is proved in \cite{Piccinini} (see also \cite{Campanato} or \cite{Kufner}), it is useful to recall the statement and the technique used in the proof, because will inspire us to techniques contained therein, for subsequent results.

\begin{prp}\label{prp1}%(\cite{Piccinini}).
\hspace*{-0.6em}\textbf{.}
%{\em
If $1<q<p<\infty, $ $0<\lambda<\mu<n, $ $q=\frac{(n\,-\,\mu)\,p}{(n\,-\,\lambda)}. $
The following embedding is true
\begin{equation}\label{4_5}
\qquad\qquad\qquad \qquad\qquad\qquad \qquad L ^{p,\lambda}(\Omega)\subset L^{q,\mu}(\Omega).
\end{equation}
\end{prp}

Proof.
%It is proved in the Piccinini's thesis \cite{Piccinini}, but it is useful to recall the technique, used later %by the author.
Applying H\"older inequality, we have
\begin{equation}\label{5}
\!\!\!\!\!\!\int_{\Omega \cap B_\rho(x)} |f|^q(y) dy \leq \left(\int_{\Omega \cap B_\rho(x)} |f|^{q\cdot \frac{p}{q}}(y) dy \right)^{\frac{q}{p}}\!\!\cdot |B_\rho|^{1\,-\,\frac{q}{p}}= C \left(\int_{\Omega \cap B_\rho(x)} |f|^{p}(y) dy \right)^{\frac{q}{p}} \!\!\!\!\cdot \!\!\rho^{n\cdot (1\,-\,\frac{q}{p})}\!\!=\,
\end{equation}
\begin{equation}\label{6}
 = C \rho^{n\cdot (1\,-\,\frac{q}{p})} \cdot \left(\frac{1}{\rho^\lambda}\int_{\Omega \cap B_\rho(x)} |f|^{p}(y) dy \right)^{\frac{q}{p}} \cdot \rho^{\lambda\cdot \frac{q}{p}}\,\leq\,
\end{equation}
\begin{equation}\label{7}
 \leq C \rho^{n\,-\,n\cdot \frac{q}{p}\,+\,\lambda\cdot \frac{q}{p}} \cdot \,\|f\|_{L^{p,\lambda}(\Omega)}=\,
\end{equation}
\begin{equation}\label{8}
 = C \rho^{\mu} \cdot \,\|f\|^q_{L^{p,\lambda}(\Omega)};
\end{equation}
then we obtain
\begin{equation}\label{9}
\frac{1}{\rho^\mu} \int_{\Omega \cap B_\rho(x)} |f|^{q}(y) dy  \leq  C \cdot \,\|f\|^q_{L^{p,\lambda}(\Omega)},
\end{equation}
where $$ \mu \,=\, n\,-\,n\cdot \frac{q}{p}\,+\,\lambda\cdot \frac{q}{p}$$
and, obviously, we have
 $$  \frac{n\,-\, \mu }{n\,-\, \lambda}=\frac{q}{p} $$
 and the conclusion follows.

\begin{rem}\label{rem2.2.}\hspace*{-0.6em}\textbf{.}
It is possible to extend the previous result considering $1 \leq q\leq p<\infty $ and $0 \leq \!\lambda, \mu \!<n $ such that $  \frac{n \,-\,\mu}{q}\,\geq \frac{n\,-\, \lambda}{p}. $
\end{rem}

\section{Embedding Results} %Section 3
\setcounter{equation}{0}
\setcounter{thm}{0}

\begin{thm}\label{them1}\hspace*{-0.6em}\textbf{.}
Let  $1<p<+\infty, $ $0<\lambda <n, \,$
$1<q<q_1<\infty, $ $0<\mu_1<\mu <1\, $
%or $1< \mu_1<\mu<n \,\,$
 and
$q=\frac{(%n
1\,-\,\mu)\,q_1}{(%n
1\,-\,\mu_1)\,\,}, $ we have
\begin{equation}\label{10}
L ^{q_1,\mu_1}(0,T,L^{p,\lambda}(\Omega))\subset L^{q,\mu}(0,T,L^{p,\lambda}(\Omega)).
\end{equation}
\end{thm}

Proof. Let us suppose that $f \in L ^{q_1,\mu_1}(0,T,L^{p,\lambda}(\Omega)), $ then is finite
\begin{equation}\label{11}
\!\!\!\!\!\!\!\!\!\!\!\!\!\!\!\!
\left(
\sup_{%\substack
{t_0 \in (0,T)\,\atop  \rho>0 }
}\,
\frac{1}{\rho^{\mu_1}}
\int_{(0,T) \cap (t_0-\rho, t_0+\rho)}
\left(
\sup_{%\substack
{x \in \Omega
\atop  %0<\rho< diam\, \Omega
\rho > 0}} \,\frac{1}{\rho^\lambda}
\int_{\Omega \cap B_\rho (x)} |f(y,t)|^p\,dy
\right)^{\!\!\frac{q_1}{p}}
dt\,
\right)^{\!\!\frac{1}{q_1}}. %\quad \forall t \in (t_0-\rho, t_0+\rho).
\end{equation}
Let us set $t \in (0,T) $ and apply H\"older inequality
\begin{equation}\label{12}
\!\!\!\!\!\!\!\!\!\!\!\!\!\!\!\!
%\left(
%\sup_{%\substack
%{t_0 \in (0,T)\,\\ \rho>0 }
%}\,
%\frac{1}{\rho^{\mu_1}}
\int_{(0,T) \cap (t_0-\rho, t_0+\rho)}
\left(
\sup_{%\substack
{x \in \Omega
\atop  %0<\rho< diam\, \Omega
\rho > 0}} \,\frac{1}{\rho^\lambda}
\int_{\Omega \cap B_\rho (x)} |f(y,t)|^p\,dy
\right)^{\!\!\frac{q}{p}}
dt\,
%\right)^{\!\!\frac{1}{q_1}}\!\!\!
\leq
\end{equation}
%errato qui sdi seguito scrivere  $|(0,T) \cap (t_0-\rho; t_0+\rho)|$?
\begin{equation}\label{13}
\!\!\!\!\!\!\!\!\!\!\!\!\!\!\!\!
\leq
\left(
%\sup_{%\substack
%{t_0, t \in (0,T)\,\\ \rho>0 ?}
%}\,
%\frac{1}{\rho^{\mu_1}}
\int_{(0,T) \cap (t_0-\rho, t_0+\rho)}
\left(
\sup_{%\substack
{x \in \Omega \atop
% 0<\rho< diam\, \Omega
\rho > 0
}} \,\frac{1}{\rho^\lambda}
\int_{\Omega \cap B_\rho (x)} |f(y,t)|^p\,dy
\right)^{\!\!\frac{q}{p}\cdot \frac{q_1}{q}}
dt\,
\right)^{\!\!\frac{q}{q_1}} \,|%B_\rho
(0,T) \cap (t_0-\rho; t_0+\rho)
|^{1\,-\,\frac{q}{q_1}}=\!\!\!
=
\end{equation}
\begin{equation}\label{14}
\!\!\!\!\!\!\!\!\!\!\!\!\!\!\!\!
=C \left(
%\sup_{%\substack
%{t_0, t \in (0,T)\,\\ \rho>0 ?}
%}\,
%\frac{1}{\rho^{\mu_1}}
\int_{(0,T) \cap (t_0-\rho, t_0+\rho)}
\left(
\sup_{%\substack
{x \in \Omega
\atop  %0<\rho< diam\, \Omega
\rho > 0
}} \,\frac{1}{\rho^\lambda}
\int_{\Omega \cap B_\rho (x)} |f(y,t)|^p\,dy
\right)^{\!\!\frac{q}{p}\cdot \frac{q_1}{q}}
dt\,
\right)^{\!\!\frac{q}{q_1}}\!\!\!\cdot \rho^{%n\cdot
(1\,-\,\frac{q}{q_1})}=\!\!\!
=
\end{equation}
\begin{equation}\label{15}
\!\!\!\!\!\!\!\!\!\!\!\!\!\!\!\!
=C \rho^{%n\cdot
(1\,-\,\frac{q}{q_1})}
\left(
%\sup_{%\substack
%{t_0, t \in (0,T)\,\\ \rho>0 ?}
%}\,
\frac{1}{\rho^{\mu_1}}
\int_{(0,T) \cap (t_0-\rho, t_0+\rho)}
\left(
\sup_{%\substack
{x \in \Omega \atop
% 0<\rho< diam\, \Omega
\rho > 0
}} \,\frac{1}{\rho^\lambda}
\int_{\Omega \cap B_\rho (x)} |f(y,t)|^p\,dy
\right)^{\!\!\frac{q_1}{p}}
dt\,
\right)^{\!\!\frac{q}{q_1}} \!\!\!\cdot \rho^{\mu_1\cdot \frac{q}{q_1}}=\!\!\!
\end{equation}
\begin{equation}\label{16}
\!\!\!\!\!\!\!\!\!\!\!\!\!\!\!\!
=C \|f\|^q_{L^{q_1,\mu_1(0,T,L^{p,\lambda}(\Omega))}}
\cdot
\rho^{1\,-\,\frac{q}{q_1}+ \mu_1\cdot\frac{q}{q_1}}.
\end{equation}
Let
\begin{equation}\label{17}
\mu\,=\, 1\,-\,\frac{q}{q_1}\,+\,\mu_1\cdot \frac{q}{q_1}= 1\,-\,(1\,-\,\mu_1)\frac{q}{q_1},
%
%=C \|f\|^q_{L^{q_1,\mu_1(0,T,L^{p,\lambda}(\Omega))}}
%\cdot
%\rho^{1\,-\,\frac{q}{q_1}+ \mu_1\cdot\frac{q}{q_1}}
\end{equation}
\begin{equation}\label{18}
\frac{1 \,-\,\mu}{1\,-\,\mu_1}\,=\, \frac{q}{q_1};
\end{equation}
it follows, as request, that
\begin{equation}\label{19}
q\,=\, \frac{(1 -\,\mu)q_1}{1\,-\,\mu_1\,\,},
\end{equation}
and the proof is complete.

\begin{rem}\label{rem2}\hspace*{-0.6em}\textbf{.}
It is possible to extend the previous result considering $1< q\leq q_1< \infty, $ $0< \mu_1\leq\mu < 1 $
 or $1<\mu_1\leq\mu<n \,\, $
 and
\begin{equation}\label{19}
\qquad \qquad \qquad \qquad\qquad \frac{1 \,-\,\mu}{q}\,\geq \, \frac{1\,-\,\mu_1}{q_1}.
\end{equation}
\end{rem}

\begin{thm}\label{them2}\hspace*{-0.6em}\textbf{.}
Let $1<q<p<\infty, \,$ $0<\lambda<\mu<n, \,\,$
$q=\frac{(n\,-\,\mu)\,p}{(n\,-\,\lambda)\,\,}, $
$1<q_2<q_1<\infty, $
$0<\mu_1<\mu_2 <1\, $
or $1< \mu_1<\mu_2 <n\,\,$ and
$q_2\,=\, \frac{(1 -\,\mu_2)q_1}{(1\,-\,\mu_1)\,\,}, $
%$q_1=\frac{(%n
%1\,-\,\mu_1)\,q_2}{(%n
%1\,-\,\mu_2)\,\,}, $
we have
\begin{equation}\label{10}
L ^{q_1,\mu_1}(0,T,L^{p,\lambda}(\Omega))\subset L^{q_2,\mu_2}(0,T,L^{q,\mu}(\Omega)).
\end{equation}
\end{thm}

Proof. %From \ref{prp1}, i
Let us set $t \in (0,T).\, $
If $1<q<p<\infty, $ $0<\lambda<\mu<n\, $ and $\, q=\frac{(n\,-\,\mu)\,p}{(n\,-\,\lambda)},\,\, $
we have, from Proposition \ref{prp1},

\begin{equation}\label{11}
\frac{1}{\rho^\mu}\int_{\Omega \cap B_\rho(x)} |f|^q(y,t) dy \leq C \left(
\sup_{%\substack
{x \in \Omega
\atop  %0<\rho< diam\, \Omega
\rho > 0
}} \,
\frac{1}{\rho^\lambda}\int_{\Omega \cap B_\rho(x)} |f|^{q\cdot \frac{p}{q}}(y,t) dy \right)^{\frac{q}{p}}.
%\cdot |B_\rho|^{1\,-\,\frac{q}{p}}\,=\, = C \cdot \left(\int_\Omega |f|^{p}(x) dx \right)^{\frac{q}{p}} \cdot %\rho^{n\cdot (1\,-\,\frac{q}{p})}\,=\,
\end{equation}
Let us fix  $t_0 \in (0,T), $ then, integrating in $(0,T) \cap (t_0-\rho, t_0+\rho), $ we have
\begin{equation}\label{12}
\int_{(0,T) \cap (t_0-\rho, t_0+\rho)}
\left(
\sup_{%\substack
{x \in \Omega
\atop  %0<\rho< diam\, \Omega
\rho >0 }} \,
\frac{1}{\rho^\mu}\int_{\Omega \cap B_\rho(x)} |f|^q(y,t) dy
\right)^{\frac{1}{q}\cdot q_2}\,
dt
\leq
\end{equation}
\begin{equation}\label{13}
\leq
C
\int_{(0,T) \cap (t_0-\rho, t_0+\rho)}
\left(
\sup_{%\substack
{x \in \Omega
\atop  %0<\rho< diam\, \Omega
\rho > 0}} \,
\frac{1}{\rho^\lambda}\int_{\Omega \cap B_\rho(x)} |f|^{q\cdot \frac{p}{q}}(y,t) dy \right)^{\frac{1}{p}\cdot q_2}
dt
\leq
\end{equation}
applying H\"older inequality, we have %errato con $|B_\rho|!$
\begin{equation}\label{14}
\leq
C
\left(
\int_{(0,T) \cap (t_0-\rho, t_0+\rho)}
\left(
\sup_{%\substack
{x \in \Omega
\atop  %0<\rho< diam\, \Omega
\rho > 0}}
\,
\frac{1}{\rho^\lambda}\int_{\Omega \cap B_\rho(x)} |f|^{p}(y,t) dy \right)^{\frac{q_2}{p}\cdot \frac{q_1}{q_2}}
dt
\right)^{\frac{q_2}{q_1}}
\!\!\!\!
\cdot %|B_\rho|
\rho^{1\,-\,\frac{q_2}{q_1}}\,=\,
%C \cdot \left(\int_\Omega |f|^{p}(x) dx \right)^{\frac{q}{p}} \cdot \rho^{n\cdot (1\,-\,\frac{q}{p})}\,=\,
\end{equation}
\begin{equation}\label{15}
= C
\left(
\frac{1}{\rho^{\mu_1}}
\int_{(0,T) \cap (t_0-\rho, t_0+\rho)}
\left(
\sup_{%\substack
{x \in \Omega
\atop  %0<\rho< diam\, \Omega
\rho > 0}}
\,
\frac{1}{\rho^\lambda}\int_{\Omega \cap B_\rho(x)} |f|^{p}(y,t) dy \right)^{\frac{q_1}{p}}
dt
\right)^{\frac{q_2}{q_1}}
\!\!\!\!
\cdot \rho^{1\,-\,\frac{q_2}{q_1}\,+\,\mu_1\cdot\frac{q_2}{q_1}}\,=\,
\end{equation}
\begin{equation}\label{16}
= C \|f\|^{q_2}_{L^{q_1,\mu_1}(0,T,L^{p,\lambda}(\Omega))} \cdot \rho^{\mu_2}
\end{equation}
where
\begin{equation}\label{17}
\mu_2\,=\,1\,-\,(1\,-\,\mu_1)\cdot \frac{q_2}{q_1},
\end{equation}
then
\begin{equation}\label{18}
\frac{1 \,-\,\mu_2}{1\,-\,\mu_1}\,=\, \frac{q_2}{q_1};
\end{equation}
it follows
\begin{equation}\label{19}
q_2\,=\, \frac{(1 -\,\mu_2)q_1}{(1\,-\,\mu_1)\,\,}.
\end{equation}
Then, we obtain
\begin{equation}\label{20}
\!\!\!\!\!\!\!\!\!\!\left(
\frac{1}{\rho^{\mu_2}}
\int_{(0,T) \cap (t_0-\rho, t_0+\rho)}
\left(
\sup_{%\substack
{x \in \Omega
\atop  %0<\rho< diam\, \Omega
\rho > 0}}
\,
\frac{1}{\rho^\mu}\int_{\Omega \cap B_\rho(x)} |f|^{q}(y,t) dy \right)^{\frac{q_2}{q}}
\!\!\!dt
\right)^{\frac{1}{q_2}}
\!\!\!\!\!\!
\leq C
\|f\|_{L^{q_1,\mu_1}(0,T,L^{p,\lambda}(\Omega))}
\end{equation}

and, finally
\begin{equation}\label{21}
\|f\|_{L^{q_2,\mu_2}(0,T,L^{q,\mu}(\Omega))}
\leq \,C\,\,
\|f\|_{L^{q_1,\mu_1}(0,T,L^{p,\lambda}(\Omega))}.
\end{equation}

\begin{rem}\label{rem3}\hspace*{-0.6em}\textbf{.}
It is possible to extend the previous result considering
$1<q\leq p<\infty, %p\geq q > 1
 \,$ $0<\lambda \leq \mu<n, \,\,$
$1 < q_2\leq q_1<\infty, $ $0<\mu_1\leq \mu_2 < 1 $ or
$1<\mu_2\leq \mu_1<n%\mu_1\geq \mu_2 > 1
\,\, $ and
\begin{equation}\label{22}
\frac{n \,-\,\mu}{q}\,\geq \, \frac{n\,-\,\lambda}{p}
;\,\,\qquad \frac{1 \,-\,\mu_2}{q_2}\,\geq \, \frac{1\,-\,\mu_1}{q_1}.
\end{equation}
\end{rem}

\section{%Applications of
%Mixed Morrey Spaces
Main Results
%Estimate of some integral operators in Mixed Morrey Spaces
} %Section 4
\subsection{Estimate of some integral operators } %Section 5
\setcounter{equation}{0}
\setcounter{thm}{0}
\setcounter{prp}{0}

Let $f\,\in \, L^1_{loc}(\R^n)\, $ and recall the following Hardy-Littlewood maximal function %$M\,f\,  $
\begin{equation}\label{23}
M\,f\,(x)=\,
\sup_{%\substack
{%x \in \Omega
\\ \rho\,>\,0 }} \,\frac{1}{|B_\rho(x)|}
\int_{%\Omega \cap
B_\rho (x)} |f(y)|\,dy
\end{equation}
where $B_\rho(x) \,$ is a ball centered at $x$ and with radius $\rho. $

%\begin{prp}%\label{prp1}
\begin{prp4.}
%\begin{prp}\label{prp2}\hspace*{-0.6em}\textbf{.}(see \cite{CF} Theorem 1).$\,$
$\!\!\!$.
Let $1<p<+\infty, $ $ 0 <\lambda<n.\,$ Then
\begin{equation}\label{24}
\| M\,f\|_{L^{p,\lambda}(\R^n%\Omega
)}
\leq \,C\,\,
\|f\|_{L^{p,\lambda}(\R^n%\Omega
)}
\end{equation}
where $C$ is independent of $f. $
\end{prp4.}
%\end{prp}
%\end{thm}

Let us now extend the previous result as follows.

\begin{thm}\label{thm4-1}\hspace*{-0.6em}\textbf{.}
Let  $1<p<+\infty, $ $0<\lambda <n, \,$ $1<q'<+\infty, $ $0< \mu<1 \,$ or $1<\mu<n \,$ and
$ f\in{L^{q',\mu}(0,T,L^{p,\lambda}(\R^n))}. $ Then,
\begin{equation}\label{25}
\| M\,f\|_{L^{q',\mu}(0,T,L^{p,\lambda}(\R^n))}
\leq \,C\,\,
\|f\|_{L^{q',\mu}(0,T,L^{p,\lambda}(\R^n))}.
\end{equation}
\end{thm}

Proof. %Let us consider  $x \in \R^n,\, $ w
%We recall that
Let $t \in (0,T).\,$
From \cite{FeStein} (Lemma 1, pg.111), we have %for $\forall t \in (0,T), $

%(corretto $\chi(y) $ o deve essere $\chi(y,t) ?$ Vorrei $\chi(y). $)
\begin{equation}\label{25_5}
\int_{ \R^n} |M\,f(y,t)|^p\chi(y)\,dy\,\leq\,c
\int_{ \R^n} |\,f(y,t)|^p(M\,\chi)(y)\,dy\,
\end{equation}
for any function $f $ and $\chi $ %$\chi \geq 0. $
%and $t_0,t \in (0,T)\,$  and,
the characteristic function of a ball $B_\rho(x)\subset \R^n, $
being the constant $c$  independent of $f.\,$ Then

%Using the method applied in \cite{CF}, we have
\begin{equation}\label{25_6}
\!\!\!\!\!\!\!%\frac{1}{\rho^\lambda}
 \int_{ B_\rho (x)} \!\!\!|M\,f(y,t)|^p\,dy\leq
%\frac{1}{(2\rho)^\lambda}
\int_{ B_{2\,\rho} (x)} \!\!\!|f(y,t)|^p (M\,\chi(y))dy +\sum_{k=1}^{+\infty}
%\frac{1}{(2^{k+1}\rho)^\lambda}
\int_{ B_{2^{k+1}\,\rho \backslash B_{2^{k}\rho} (x)}} \!\!\!\!\!\!\!\!\!\!\!\!\!\!|f(y,t)|^p(M\,\chi(y))\,dy
\end{equation}
%then
it follows
%\begin{equation}\label{26}
$$
\!\!\!\!\!\!\!\frac{1}{\rho^\lambda} \int_{ B_\rho (x)} |M\,f(y,t)|^p\,dy\,\leq\,
$$
%\end{equation}
\begin{equation}\label{26_5}
\leq C \frac{1}{(2\rho)^\lambda} \int_{ B_{2\,\rho} (x)} |f(y,t)|^p (M\,\chi(y))\,dy \,+C\sum_{k=1}^{+\infty}
\frac{1}{(2^{k+1}\rho)^\lambda} \int_{ B_{2^{k+1}\,\rho} (x)} \!\!\!\!\!\!|f(y,t)|^p(M\,\chi(y))\,dy,
\end{equation}

using the method applied in \cite{CF} and considering the supremum for $x \in \R^n $ and $\rho>0. $
Let us fix  $ t_0 \in (0,T), $ then,
elevating to $\frac{q'}{p}, $  integrating in $(0,T) \cap (t_0-\rho, t_0+\rho) $ and multiplying for $\rho^{-\mu}, $ we obtain
\begin{equation}\label{27}
\frac{1}{\rho^\mu} \int_{(0,T) \cap (t_0-\rho, t_0+\rho)}
\left(
\sup_{%\substack
{x \in \R^n \atop
 \rho\,>\,0 }}
\frac{1}{\rho^\lambda} \int_{ B_\rho (x)} |M\,f(y,t)|^p\,dy\,
\right)^{\frac{q'}{p}}
dt
\leq\,
\end{equation}
\begin{equation}\label{28}
\!\!\!\!\leq C \frac{1}{\rho^\mu}
\int_{(0,T) \cap (t_0-\rho, t_0+\rho)}
\left(
\sup_{%\substack
{x \in \R^n
\atop  \rho\,>\,0 }}
\frac{1}{\rho^\lambda}
\int_{ B_{\rho} (x)} |f(y,t)|^p\,dy \,
\right)^{\frac{q'}{p}}
dt
\end{equation}
taking the supremum, in both sides, for  $t_0 \in (0,T)$ and $\rho \,>\,0, $ we obtain
\begin{equation}\label{29}
\left[
\sup_{%\substack
{t_0 \in (0,T)\,\atop  \rho>0 }}\,
\frac{1}{\rho^{\mu}}
\int_{(0,T) \cap (t_0-\rho, t_0+\rho)}
%\frac{1}{\rho^\lambda} \int_{(0,T) \cap (t_0-\rho, t_0+\rho)}
\left(
\sup_{%\substack
{x \in \R^n \atop
 \rho\,>\,0 }}
\frac{1}{\rho^\lambda} \int_{ B_\rho (x)} |M\,f(y,t)|^p\,dy\,
\right)^{\frac{q'}{p}}
dt
\right]^{\frac{1}{q'}}\leq
\end{equation}

\begin{equation}\label{30}
\leq \left[
\sup_{%\substack
{t_0, t \in (0,T)\,\atop  \rho>0 }}\,
\frac{1}{\rho^{\mu}}
\int_{(0,T) \cap (t_0-\rho, t_0+\rho)}
%\frac{1}{\rho^\lambda} \int_{(0,T) \cap (t_0-\rho, t_0+\rho)}
\left(
\sup_{%\substack
{x \in \R^n
\atop  \rho\,>\,0 }}
\frac{1}{\rho^\lambda} \int_{ B_\rho (x)} |\,f(y,t)|^p\,dy\,
\right)^{\frac{q'}{p}}
dt
\right]^{\frac{1}{q'}}
\end{equation}
or, equivalently
\begin{equation}\label{31}
\| M\,f\|_{L^{q',\mu}(0,T,L^{p,\lambda}(\R^n))}
\leq \,C\,\,
\|f\|_{L^{q',\mu}(0,T,L^{p,\lambda}(\R^n))}.
\end{equation}

As application of this result we prove some estimates of the Riesz potential in $L^{q,\mu}(0,T,L^{p,\lambda}(\R^n)) $ spaces.

Let us set $t \in (0,T) $ and %recall
consider, for $0<\alpha<n, $ the fractional integral operator of order $\alpha, $

%PUO' ESSERE $f(y,t) $ ?

\begin{equation}\label{32}
I_\alpha f(x,t
)\,=\, \int_{\R^n} \frac{f(y,t)}{|x
\,-\,y|^{n\,-\,\alpha}} dy, \quad a.e.\, in \,\R^n . %forall t \in (0,T).
\end{equation}

\begin{thm}\label{thm4-2}\hspace*{-0.6em}\textbf{.}
Let  $0<\alpha<n, $ $1<p<\frac{n}{\alpha, }, $ $0<\lambda <n\,-\,\alpha\,p, \,$ $\frac{1}{q}=\frac{1}{p}\,-\,\frac{\alpha}{n-\lambda} $
$1<q'<+\infty, $ $0< \mu'<1 \,$ %(or $1<\mu'<n$ ?)
 and
$ f\in{L^{q',\mu'}(0,T,L^{p,\lambda}(\R^n))}. $ Then,
\begin{equation}\label{33}
\| I_\alpha\,f\|_{L^{q',\mu'}(0,T,L^{q,\lambda}(\R^n))}
\leq \,C\,\,
\|f\|_{L^{q',\mu'}(0,T,L^{p,\lambda}(\R^n))}.
\end{equation}
\end{thm}

Proof. Let us fix $x \in \R^n, $ $ t_0 \in (0,T) $  and $ f \in L^{q',\mu'}(0,T,L^{p,\lambda}(\R^n)). $ Then, set $t \in (0,T), $

%PUO' ESSERE $f(y,t) $ ?

\begin{equation}\label{34}
(I_\alpha f)(x,t
)= \int_{|x
-y| \leq \epsilon} \frac{f(y,t)}{|x
\,-\,y|^{n\,-\,\alpha}} dy+
\int_{|x
-y| > \epsilon} \frac{f(y,t)}{|x
\,-\,y|^{n\,-\,\alpha}} dy
= I_1\,+\,I_2, %\quad  \forall t \in (0,T),
\end{equation}
estimating separately each integral $I_1$ and $I_2, $  as in \cite{Adams} (Theorem 3.1) or \cite{CF} (Theorem 2), we  obtain
%for every $\forall t \in (0,T), $

%PUO' ESSERE $(y,t),$ AL POSTO DI $(x)$, AL PRIMO  AL SECONDO MEMBRO?
\begin{equation}\label{35}
|I_\alpha f|(x,t
) \leq  C
(Mf)^{\frac{n-\lambda-\alpha p}{n-\lambda}}(x
) \cdot
\left(
\sup_{%\substack
{x \in \R^n \atop  \rho\,>\,0 }}
\frac{1}{\rho^\lambda} \int_{B_\rho (x)} |f(y,t)|^p dy
\right)^{\frac{1}{p}\cdot \frac{\alpha p}{n-\lambda}}\!\!\!\!\!\!\!\!\!\!\!\!\!\!,
\end{equation}
recalling that $\frac{n-\lambda-\alpha p}{n-\lambda} = \frac{p}{q}$, elevating to the power $q, $   integrating in $B_\rho(x) $ and multiplying to $\rho^{-\lambda }, $  we have
\begin{equation}\label{36}
\frac{1}{\rho^\lambda} \int_{B_\rho(x)} |I_\alpha\,f\,(y,t)|^q%{\frac{q}{p}}
 dy \leq
\frac{1}{\rho^\lambda} \int_{B_\rho(x)}
(Mf)^{p}(y,t) \cdot
\left(
\sup_{%\substack
{x \in \R^n \atop  \rho\,>\,0 }}
\frac{1}{\rho^\lambda} \int_{B_\rho (x)} |f(y,t)|^p dy
\right)^{\frac{\alpha q}{n-\lambda}} \!\!\!\!\!\!\!\!dy
\leq
\end{equation}
applying Theorem \ref{thm4-1}, and observing that $\frac{\alpha\,q}{n-\lambda} +1 =\frac{q}{p}, $
%(si puo' uscire il primo integrale fuori dal secondo ? )
\begin{equation}\label{37}
\leq C
\left(
\sup_{%\substack
{x \in \R^n \atop  \rho\,>\,0 }}
\frac{1}{\rho^\lambda} \int_{B_\rho (x)} |f(y,t)|^p dy
\right)^{\frac{\alpha q}{n-\lambda}} \!\!\!\!\!\!
\cdot
\left(
\sup_{%\substack
{x \in \R^n \atop  \rho\,>\,0 }}
\frac{1}{\rho^\lambda} \int_{B_\rho (x)} |f(y,t)|^p dy
\right) =
\end{equation}
\begin{equation}\label{38}
=\,C\,
\left(
\sup_{%\substack
{x \in \R^n \atop  \rho\,>\,0 }}
\frac{1}{\rho^\lambda} \int_{B_\rho (x)} |f(y,t)|^p dy
\right)^{\frac{q}{p}},
\end{equation}
then
\begin{equation}\label{39}
\frac{1}{\rho^\lambda} \int_{B_\rho(x)} |I_\alpha\,f\,(y,t)|^q%{\frac{q}{p}}
 dy \leq C
 \left(
\sup_{%\substack
{x \in \R^n \atop  \rho\,>\,0 }}
\frac{1}{\rho^\lambda} \int_{B_\rho (x)} |f(y,t)|^p dy
\right)^{\frac{q}{p}}\!\!\!,
\end{equation}
considering the supremum for $x \in \R^n$ and $\rho > 0 $ and elevating both member to $\frac{1}{q} $
\begin{equation}\label{40}
\!\!\!\!\!\left(
\sup_{%\substack
{x \in \R^n \atop  \rho\,>\,0 }}
\frac{1}{\rho^\lambda} \int_{B_\rho(x)} |I_\alpha\,f\,(y,t)|^q%{\frac{q}{p}}
 dy
 \right)^{\frac{1}{q}}
 \leq C
 \left(
\sup_{%\substack
{x \in \R^n \atop  \rho\,>\,0 }}
\frac{1}{\rho^\lambda} \int_{B_\rho (x)} |f(y,t)|^p dy
\right)^{\frac{1}{p}}\!\!\!.
\end{equation}
Now, elevating to $q',$ integrating in $(0,T)\cap (t_0-\rho, t_0+\rho) $ and multiplying to $\rho^{-\mu}, $ we have
\begin{equation}\label{41}
%\left[
%\sup_{%\substack
%{t_0, t \in (0,T)\,\\ \rho>0 ?}}\,
\frac{1}{\rho^{\mu}}
\int_{(0,T) \cap (t_0-\rho, t_0+\rho)}
\left(
\sup_{%\substack
{x \in \R^n \atop  \rho\,>\,0 }}
\frac{1}{\rho^\lambda} \int_{B_\rho(x)} |I_\alpha\,f\,(y,t)|^q
 dy
 \right)^{\frac{q'}{q}}
dt
\leq
\end{equation}

\begin{equation}\label{42}
\leq
C
\frac{1}{\rho^{\mu}}
\int_{(0,T) \cap (t_0-\rho, t_0+\rho)}
\left(
\sup_{%\substack
{x \in \R^n \atop  \rho\,>\,0 }}
\frac{1}{\rho^\lambda} \int_{B_\rho (x)} |f(y,t)|^p dy
\right)^{\frac{q'}{p}}\!\!\!dt,
\end{equation}

taking the supremum for $t_0\in (0,T),  \rho>0, $ we have
\begin{equation}\label{43}
\sup_{%\substack
{t_0 \in (0,T)\atop  \rho>0 }}\,
\frac{1}{\rho^{\mu}}
\int_{(0,T) \cap (t_0-\rho, t_0+\rho)}
\left(
\sup_{%\substack
{x \in \R^n \atop  \rho\,>\,0 }}
\frac{1}{\rho^\lambda} \int_{B_\rho(x)} |I_\alpha\,f\,(y,t)|^q
 dy
 \right)^{\frac{q'}{q}}
dt
\leq
\end{equation}
\begin{equation}\label{44}
\leq
C
\sup_{%\substack
{t_0 \in (0,T)\atop  \rho>0 }}\,
\frac{1}{\rho^{\mu}}
\int_{(0,T) \cap (t_0-\rho, t_0+\rho)}
\left(
\sup_{%\substack
{x \in \R^n \atop  \rho\,>\,0 }}
\frac{1}{\rho^\lambda} \int_{B_\rho (x)} |f(y,t)|^p dy
\right)^{\frac{q'}{p}}\!\!\!dt.
\end{equation}
Finally, elevating to $\frac{1}{q'} $ we have
\begin{equation}\label{45}
\| I_\alpha\,f\|_{L^{q',\mu'}(0,T,L^{q,\lambda}(\R^n))}
\leq \,C\,\,
\|f\|_{L^{q',\mu'}(0,T,L^{p,\lambda}(\R^n))}.
\end{equation}

\begin{crlr}\label{thm4-3}\hspace*{-0.6em}\textbf{.}
Let  $0<\alpha<n, $ $1<p<\frac{n}{\alpha, }, $ $0<\lambda <n\,-\,\alpha\,p. \,$
%(non $0<\lambda <n$ come nel %Teorema precedente e nel "Corollary" in [CF], ok? ).
%
%$\frac{1}{q}=\frac{1}{p}\,-\,\frac{\alpha}{n-\lambda} $ (al denominatore non $n$ come nel % Teorema %precedente e nel "Corollary" in [CF] ?)
%$\mu = \frac{n \lambda}{(n - \alpha)}$ (al den. non $(n-\alpha p), $ vero ?) .
Let us also set $ 1<q<p $ such that $\frac{1}{q}=\frac{1}{p}\,-\,\frac{\alpha}{n%-\lambda
},\, $
%(al den. con $n$ come nel "Corollary" in [CF] ?),
$\lambda<\mu<n\, $ %$q=\frac{(n\,-\,\mu)\,p}{(n\,-\,\lambda)}\, $
such that $\mu = \frac{n \lambda}{(n - \alpha p)},\,$
%(al den. con $(n-\alpha p), $  non con  " $n-\alpha$" come in "Corollary  [CF]", vero ?)
%
 $ \,1<q'<+\infty, $ $0< \mu'<1 \,$ %(or $1<\mu'<n $ ?)
 and
 $ f\in{L^{q',\mu'}(0,T,L^{p,\lambda}(\R^n))}. $ Then,
\begin{equation}\label{46}
\| I_\alpha\,f\|_{L^{q',\mu'}(0,T,L^{q,\mu}(\R^n))}
\leq \,C\,\,
\|f\|_{L^{q',\mu'}(0,T,L^{p,\lambda}(\R^n))}.
\end{equation}
where $C\,$ is independent of $f.$
\end{crlr}

Proof. Let us fix $x \in \R^n $ and $t%_0
 \in (0,T) $. % (da confermare).
% $ f \in L^{q',\mu'}(0,T,L^{p,\lambda}(\R^n)). $ Just observe,
%from Proposition \ref{prp1}%(\cite{Piccinini})
%, that, if $1<q<p<\infty, $ $0<\lambda<\mu<n, $ $q=\frac{(n\,-\,\mu)\,p}{(n\,-\,\lambda)}, $
%the following embedding is true
%\begin{equation}\label{4_5}
%  L ^{p,\lambda}(\R^n)\subset L^{q,\mu}(\R^n).
%\end{equation}
%
%Then, f
From %Theorem 2
 Corollary in \cite{CF}, %$\forall t \in (0,T), $
  we have
%PUO' ESSERE $(y,t),$ AL POSTO DI $(y)$, AL PRIMO  AL SECONDO MEMBRO? SI, anzi devo.

\begin{equation}\label{47}
\left(
\sup_{%\substack
{x \in \R^n \atop  \rho\,>\,0 }}
\frac{1}{\rho^\mu} \int_{B_\rho(x)} |I_\alpha\,f\,(y,t)|^q
 dy
 \right)^{\frac{1}{q}}
\leq \,C\,
\left(
%C
\,
\sup_{%\substack
{x \in \R^n \atop  \rho\,>\,0 }}
\frac{1}{\rho^\lambda} \int_{B_\rho (x)} |f(y,t)|^p dy
\right)^{\frac{1}{p}},
\end{equation}
elevating to $q', $
 fixing $t_0 \in (0,T),\, $
integrating in $(0,T) \cap (t_0-\rho, t_0+\rho) $  and multiplying for $\rho^{-\mu'}, $
\begin{equation}\label{48}
\left(
\sup_{%\substack
{t_0 \in (0,T)\atop  \rho>0 }}\,
\frac{1}{\rho^{\mu'}}
\int_{(0,T) \cap (t_0-\rho, t_0+\rho)}
\left(
\sup_{%\substack
{x \in \R^n \atop  \rho\,>\,0 }}
\frac{1}{\rho^\lambda} \int_{B_\rho(x)} |I_\alpha\,f\,(y,t)|^q
 dy
 \right)^{\frac{q'}{q}}
dt
\right)^{\frac{1}{q'}}
\leq
\end{equation}
\begin{equation}\label{49}
\leq
C
\left(
\sup_{%\substack
{t_0 \in (0,T)\atop  \rho>0 }}\,
\frac{1}{\rho^{\mu'}}
\int_{(0,T) \cap (t_0-\rho, t_0+\rho)}
\left(
\sup_{%\substack
{x \in \R^n \atop  \rho\,>\,0 }}
\frac{1}{\rho^\lambda} \int_{B_\rho (x)} |f(y,t)|^p dy
\right)^{\frac{q'}{p}}\!\!\!dt
\right)^{\frac{1}{q'}}
\end{equation}
that is
\begin{equation}\label{50}
\| I_\alpha\,f\|_{L^{q',\mu'}(0,T,L^{q,\mu}(\R^n))}
\leq \,C\,\,
\|f\|_{L^{q',\mu'}(0,T,L^{p,\lambda}(\R^n))}.
\end{equation}

One more application of the technique used in the proof of Theorem \ref{thm4-1} is the following result, where we set $T$  a convolution singular integral operator $T\,=\, k * f, $ where $k $  is an usual Calder\'on-Zygmund kernel, studied by Coifman and Fefferman in \cite{CoFe}.
\begin{thm}\label{thm4-3}\hspace*{-0.6em}\textbf{.}
Let  $1<p<\infty, $ $0<\lambda <n\,$ $ 1<q'<+\infty, $ $ 0< \mu'<1 \,$ %(or $1<\mu'<n ?$)
 and
$ f\in{L^{q',\mu'}(0,T,L^{p,\lambda}(\R^n))}. $ Then,
\begin{equation}\label{33}
\| T\,f\|_{L^{q',\mu'}(0,T,L^{q,\lambda}(\R^n))}
\leq \,C\,\,
\|f\|_{L^{q',\mu'}(0,T,L^{p,\lambda}(\R^n))}.
\end{equation}
\end{thm}

Proof. Let us fix $x \in \R^n, $ $t%_0
 \in (0,T), $ %(da confermare),
 $ f \in L^{q',\mu'}(0,T,L^{p,\lambda}(\R^n))$ and  $\chi $ the characteristic function of a ball $B_\rho (x). $ Then, from  a result by Coifman and Rochberg (see \cite{CoRo} pg.251),  %and  Proposition \ref{prp2},
 $M(M\chi)^\gamma \leq c\,  (M\chi)^\gamma, $ then $(M\chi)^\gamma $ is a $A_1$ weight.

It follows, from a result contained in \cite{CoFe}, %$\forall t \in (0,T), $
that %(MA $\gamma$ C'E' OPPURE NO ??)
\begin{equation}\label{34}
\!\!\!\!\!\!\!%\frac{1}{\rho^\lambda}
 \int_{ B_\rho (x)} |T\,f(y,t)|^p\,dy\leq
%\frac{1}{(2\rho)^\lambda}
\int_{ \R^n }  |T\,f(y,t)|^p (M\,\chi(y))^\gamma dy \leq C
\int_{ \R^n } |f(y,t)|^p (M\,\chi(y))^\gamma \,dy,
\end{equation}
%and, taking in mind
estimating the last term following the lines of the proof of Theorem \ref{thm4-1},  we get the conclusion.

Before we prove the next results we need to consider two variants of the Hardy-Littlewood maximal operator, that are the Sharp Maximal function and the Fractional maximal functions (see e. g. \cite{DR1}).

\begin{Def}%Definition of Sharp Maximal function vedi [DR1- Commutators in Morreey spaces]
%Let $1<p<+ \infty, $ $0<\lambda<n, $ and let $f $ be a real measurable
%function on the open bounded set $\Omega\,\subset\, \R^n.$
%%(n \geq 1).$

Given $f\in L^1_{loc} (\R^n)\,$ let us define the following Sharp Maximal function
\begin{equation}\label{35}
f^\sharp (x)\,=\, \sup_{ B \supset \{x\} } \,\frac{1}{|B|}
\int_{B} |f(y)\,-\, f_B|\,dy,
\end{equation}
for a.e. $x \in \R^n, $ where $B $ is a generic ball in $\R^n. $
\end{Def}

\begin{Def}%Definition of Sharp Maximal function vedi [DR1- Commutators in Morreey spaces]
%Let $1<p<+ \infty, $ $0<\lambda<n, $ and let $f $ be a real measurable
%function on the open bounded set $\Omega\,\subset\, \R^n.$
%%(n \geq 1).$

%Given
Set $t \in (0,T), \,$ $f\in L^1_{loc} (\R^n)\,$ and $0<\eta<1. $ Let us define the Fractional Maximal function
\begin{equation}\label{36}
(M_\eta f) (x)\,=\, \sup_{ B \supset \{x\} } \,\frac{1}{|B|^{1\,-\,\eta}}
\int_{B} |f(y,t)\,-\, f_B|\,dy,
\end{equation}
for a.e. $x \in \R^n, $ where $B $ is a generic ball in $\R^n. $
\end{Def}

The next %Proposition
 Theorem is a generalization of a well known inequality by Fefferman and Stein, see \cite{FeStein}, pg. 153.

%\begin{prp}\label{prp3}\hspace*{-0.6em}\textbf{.}
\begin{thm}\label{thmprp1}\hspace*{-0.6em}\textbf{.}
Let $1< p,q <\infty, $ $0<\lambda, \mu <n $ and $ f \in L^{q,\mu}(0,T,L^{p,\lambda}(\R^n)). $

Then, there exists a constant $C\geq 0 $ independent of $f $ such that
\begin{equation}\label{37}
\| M\,f\|_{L^{q,\mu}(0,T,L^{p,\lambda}(\R^n))}
\leq \,C\,\,
\|f^\sharp\|_{L^{q,\mu}(0,T,L^{p,\lambda}(\R^n))}.
\end{equation}
%\end{prp}
\end{thm}

Proof. Let us fix $x \in \R^n, t%_0
 \in (0,T).\,$
 Let us also consider $\rho>0, \gamma \in ]\frac{\lambda}{n}; 1[,\,$ $\chi=\chi_{B_\rho(x)},  $ $\forall x \in \R^n $ the characteristic function of a ball $B_\rho(x). $ We know
 %(see  \cite{CoRo} pg.251),
 that $(M\chi)^\gamma \in A_1 $ and, from \cite{GC} pg. 410, %$\forall t \in (0,T), $
  we have

%(HO INSERITO NEL RIGO SUCCESSIVO  LA VARIABILE "$t$", E' CORRETTO %SE NON METTO ANCHE "dt" ??)
\begin{equation}\label{38}
\int_{ \R^n }  (M\,f%(y,t)
)^p (y,t) \omega (y%,t
) dy \leq C
\int_{ \R^n } |f^\sharp(y,t)|^p  \omega (y%,t
) \,dy,\quad \forall \omega \in A_\infty, \forall f \in L^p_\omega (\R^n)
\end{equation}
where $L^p_\omega(\R^n) $ is the $L^p$ space with respect to the measure $d\mu =\omega\,dx. $ We can use this inequality   because $f \in L^{p,\lambda}(\R^n) $ implies  $f \in L^p_{(M\chi)^\gamma} (\R^n)$ (see the calculation in \cite {CF} pg. 275).

Choosing $\omega(y%,t
)= (M \chi)^\gamma (y), $ we have, from \cite{DR1} pg.327, %???????????????????????????????????????????????????????

\begin{equation}\label{39}
\int_{ B_\rho(x) }  (M\,f%(y,t)
)^p (y,t) dy \leq
\int_{ \R^n }  (M\,f)^p (y,t)  (M\,\chi)^\gamma (y) dy \leq
\end{equation}
\begin{equation}\label{39,5}
\leq C \cdot \int_{ \R^n} |f^\sharp(y,t)|^p (M\,\chi)^\gamma (y)  \,dy \leq
C \rho^\lambda \,
\sup_{%\substack
{x \in \R^n \atop  \rho\,>\,0 }}
\frac{1}{\rho^\lambda}
\int_{ B_\rho(x) } |f^\sharp(y,t)|^p  \,dy,\quad  \forall f \in L^p_\omega (\R^n)
\end{equation}
then
\begin{equation}\label{40}
\frac{1}{\rho^\lambda}
\int_{ B_\rho(x) }  (M\,f)^p (y,t) dy \leq C \,
\sup_{%\substack
{x \in \R^n \atop  \rho\,>\,0 }}
\frac{1}{\rho^\lambda}
\int_{ B_\rho(x) } |f^\sharp(y,t)|^p  \,dy,\quad  \forall f \in L^p_\omega (\R^n)
\end{equation}
and, taking the supremum for $x \in \R^n $ and $ \rho\,>\,0 $ we have
\begin{equation}\label{41}
\sup_{%\substack
{x \in \R^n \atop  \rho\,>\,0 }}
\frac{1}{\rho^\lambda}
\int_{ B_\rho(x) }  (M\,f)^p (y,t) dy \leq C
\sup_{%\substack
{x \in \R^n \atop  \rho\,>\,0 }}
\frac{1}{\rho^\lambda}
\int_{ B_\rho(x) } |f^\sharp(y,t)|^p  \,dy,\quad  %\forall f \in L^p_\omega (\R^n)
\end{equation}
set $t_0 \in (0,T), $
elevating to $\frac{q}{p}, $  integrating in $(0,T) \cap (t_0-\rho, t_0+\rho) $ and multiplying for $\rho^{-\mu}, $ we have
\begin{equation}\label{42}
\frac{1}{\rho^\mu} \int_{(0,T) \cap (t_0-\rho, t_0+\rho)}
\left(
\sup_{%\substack
{x \in \R^n \atop  \rho\,>\,0 }}
\frac{1}{\rho^\lambda}
\int_{ B_\rho(x) }  (M\,f(y,t))^p (y,t) dy
\right)^{\frac{q}{p}}
dt
\leq\,
\end{equation}
\begin{equation}\label{43}
\leq C
\frac{1}{\rho^\mu} \int_{(0,T) \cap (t_0-\rho, t_0+\rho)}
\left(
\sup_{%\substack
{x \in \R^n \atop  \rho\,>\,0 }}
\frac{1}{\rho^\lambda}
\int_{ B_\rho(x) } |f^\sharp(y,t)|^p  \,dy
\right)^{\frac{q}{p}}
dt
\end{equation}
then, we obtain
\begin{equation}\label{44}
\left(
\sup_{%\substack
{t_0 \in (0,T)\atop  \rho>0 }}\,
\frac{1}{\rho^\mu} \int_{(0,T) \cap (t_0-\rho, t_0+\rho)}
\left(
\sup_{%\substack
{x \in \R^n \atop  \rho\,>\,0 }}
\frac{1}{\rho^\lambda}
\int_{ B_\rho(x) }  (M\,f(y,t))^p (y,t) dy
\right)^{\frac{q}{p}}
dt
\right)^{\frac{1}{q}}
\leq\,
\end{equation}
\begin{equation}\label{45}
\leq C
\left(
\sup_{%\substack
{t_0 \in (0,T)\atop  \rho>0 }}\,
\frac{1}{\rho^\mu} \int_{(0,T) \cap (t_0-\rho, t_0+\rho)}
\left(
\sup_{%\substack
{x \in \R^n \atop  \rho\,>\,0 }}
\frac{1}{\rho^\lambda}
\int_{ B_\rho(x) } |f^\sharp(y,t)|^p  \,dy
\right)^{\frac{q}{p}}
dt
\right)^{\frac{1}{q}}
\end{equation}
and we get the conclusion.

%\begin{prp}\label{prp4}\hspace*{-0.6em}\textbf{.} QUII
\begin{thm}\label{thmprp2}\hspace*{-0.6em}\textbf{.}
Let $1< p, q, q_1 <\infty, $ $0<\lambda, \mu_1 <n $ and $ f \in L^{q_1,\mu_1}(0,T,L^{p,\lambda}(\R^n)). $

Then, for every $\eta \in ]0,(1-\frac{\lambda}{n})\frac{1}{p}[, $ there exists a constant $C\geq 0 $ independent of $f $ such that
\begin{equation}\label{46}
\| M_\eta\,f\|_{L^{q_1,\mu_1}(0,T,L^{q,\lambda}(\R^n))}
\leq \,C\,\,
\|f%^\sharp
\|_{L^{q_1,\mu_1}(0,T,L^{p,\lambda}(\R^n))}
\end{equation}
where
\begin{equation}\label{47}
\frac{1}{q}\,=\,\frac{1}{p}\,-\,\frac{n\,\eta}{n\,-\,\lambda}.
\end{equation}
%\end{prp} QUII
\end{thm}
Proof. Let $x \in \R^n\,$ and $t_0 \in (0,T). $

Let us fix $1<r<p$ and
\begin{equation}\label{48}
\varepsilon\,=\,\,\frac{\left(1\,-\,\frac{\lambda}{n} \right) \cdot \frac{p}{n}\,-\,\eta}{\left(1\,-\,\frac{\lambda}{n} \right) \frac{1}{p}}.
\end{equation}
Set $t \in (0,T), $ %$\forall t \in (0,T) $ 
for a generic ball $B$ of $\R^n, $ we have
\begin{equation}\label{49}
\,\frac{1}{|B|^{1\,-\,\eta}} \int_B\, |f(y,t)|dy\leq
\end{equation}
\begin{equation}\label{50}
\leq \left(
\,\frac{1}{|B|} \int_B\, |f(y,t)|^r dy \right)^{\frac{\varepsilon}{r}}
\cdot
\left( \,\frac{1}{|B|^{\frac{\lambda}{n} }} \int_B\, |f(y,t)|^p dy \right)^{\frac{(1-\varepsilon)}{p}}
\end{equation}
then

%(E' CORRETTO SCRIVERE $(y,t)$ SUBITO SOTTO IN $  \left[M(|f(y,t)|^r)  \right]^{\frac{\varepsilon}{r}} (x) $ ?)
\begin{equation}\label{51}
\,\frac{1}{|B|^{1\,-\,\eta}} \int_B\, |f(y,t)|dy\leq
\,\left[  M(|f|^r)  \right]^{\frac{\varepsilon}{r}} (y,t)
\cdot
\left(
\sup_{%\substack
{x \in \R^n \atop  \rho\,>\,0 }}
\,\frac{1}{\rho^{\lambda}} \int_B\, |f(y,t)|^p dy
\right)^{1-\varepsilon }
\end{equation}
from which it follows
\begin{equation}\label{52}
\left(  M_\eta(f)  \right)^{\frac{p}{\varepsilon}}(y,t)
\leq
\left(  M(|f|^r)  \right)^{\frac{p}{r}}(y,t)
\cdot
\|f\|_{L^{p,\lambda}(\R^n)}^{\frac{(1\,-\,\varepsilon)}{\varepsilon}\cdot p}\qquad\,\,a.\,\,e.\,\, y \in\R^n,\,t \in (0,T).
\end{equation}
Denoting by $\chi(y%,t
)\,=\,\chi_{B_\rho(x)}(y)\,$ we have

%AL I MEMBRO E' CORRETTO $\left(  M_\eta(f)  \right)^{\frac{p}{\varepsilon}}(y,t)$ OPPURE
%$ \left(  M_\eta(f(y,t))  \right)^{\frac{p}{\varepsilon}}(x)$ ?
%
%SIMILMENTE AL II MEMBRO: $\left(  M(|f|^r)  \right)^{\frac{p}{r}}(y,t)$ OPPURE
%$\left(  M(|f\left(  M(|f|^r)  \right)^{\frac{p}{r}}(y,t)|^r)  \right)^{\frac{p}{r}}\!(x)$ ?
\begin{equation}\label{53}
\int_{\R^n}
\left(  M_\eta(f)  \right)^{\frac{p}{\varepsilon}}(y,t)
\cdot \chi(y%,t
)\,dy
\leq
\|f\|_{L^{p,\lambda}(\R^n)}^{\frac{(1\,-\,\varepsilon)}{\varepsilon}\cdot p}\,\,
\int_{\R^n}
\left(  M(|f|^r)  \right)^{\frac{p}{r}}(y,t)
\cdot \chi(y%,t
)\,dy
\leq
\end{equation}
\begin{equation}\label{54}
%\int_{\R^n}
%\left(  M_\eta(f)  \right)^{\frac{p}{\varepsilon}}(y,t)
%\cdot \chi(y,t)\,dy
\qquad\qquad\qquad\qquad
\qquad\qquad\qquad\leq
\,
\|f\|_{L^{p,\lambda}(\R^n)}^{\frac{(1\,-\,\varepsilon)}{\varepsilon}\cdot p}\,\,
\int_{\R^n}
\, |f|^p\,(y,t) \,\,
\cdot %M\, \chi(y)
(M\, \chi(y))%^\gamma
\,\,dy.
\end{equation}
%(nel rigo precedente e' $M\, \chi(y)$  oppure $(M\, \chi(y))^\gamma $ come in Lemma 3 in \cite{DR1} ?)

%Estimating the last integral as in the previous Proposition %\ref{prp3}
% we get
Then, we obtain
\begin{equation}\label{55}
\int_{B_\rho(x) }
\left(  M_\eta(f)  \right)^{\frac{p}{\varepsilon}}(y,t)\,dy
\leq \,C\,
\|f\|_{L^{p,\lambda}(\R^n)}^{%\frac{(1\,-\,\varepsilon)}{\varepsilon}%\cdot p
\frac{p}{\varepsilon}}
\,\,
%\int_{\R^n}
%\left(  M(|f|^r)  \right)^{\frac{p}{r}}(y,t)
%\cdot \chi(y,t)\,dy
\cdot
\rho^\lambda.
\end{equation}

Let us observe that 
%$$
\begin{equation}\label{55,5}
\frac{p}{\varepsilon}\,=\, q
\end{equation}
%$$
indeed, %in the statement, see  from 
using (\ref{47}), we have %posed
%\begin{equation}\label{56}
%\frac{1}{q}\,=\,\frac{1}{p}\,-\,\frac{n\,\eta}{n\,-\,\lambda},
%\end{equation}
%then
%\begin{equation}\label{57}
%\frac{\varepsilon}{p}\,=\,\frac{1}{p}\,-\,\frac{n\,\eta}{n\,-\,\lambda},
%\end{equation}
%or, equivalently
%\begin{equation}\label{58}
%\varepsilon=\,1\,-\,\frac{n\,\eta\,p}{n\,-\,\lambda},
%\end{equation}
%then
\begin{equation}\label{58}
\varepsilon=\,\,\frac{n\,-\,\lambda\,-n\,\eta\,p}{n\,-\,\lambda},
\end{equation}
dividing by $n \,\cdot \,p, $ we %have
%\begin{equation}\label{59}
%\varepsilon=\,\,\frac{1\,-\,\frac{\lambda}{n}\,-\,\eta\,p}{1\,-\,\frac{\lambda}{n}},
%\end{equation}
%and also dividing by $p$
%\begin{equation}\label{60}
%\varepsilon=\,\,\frac{(1\,-\,\frac{\lambda}{n})\frac{1}{p}\,-\,\eta\,}{(1\,-\,\frac{\lambda}{n})\frac{1}{p}},
%\end{equation}
%the last %in
%equality is true because it is
deduce exactly (\ref{48}). %what we fix  at the beginning of the proof.

Then, we obtain %that $\frac{p}{\varepsilon}= q\,$ and
\begin{equation}\label{61}
\left(
\sup_{%\substack
{x \in \R^n \atop  \rho\,>\,0 }}
\frac{1}{\rho^\lambda}
\int_{B_\rho(x) }
\left(  M_\eta(f)  \right)^{%\frac{p}{\varepsilon}
q}(y,t)\,dy
\right)^{\frac{1}{q}}
\leq \,
C\,
\left(
\sup_{%\substack
{x \in \R^n \atop  \rho\,>\,0 }}
\frac{1}{\rho^\lambda}
\int_{B_\rho(x) }
  |f|^p (y,t) \,dy
\right)^{\frac{1}{p%q
}},
\end{equation}
%set $t_0 \in (0,T), $
elevating to $q_1\,$ integrating both sides in $(0,T)\cap (t_0-\rho;t_0+\rho) $ and  multiplying for $\frac{1}{\rho^{\mu_1}}, $ we have
\begin{equation}\label{62}
\frac{1}{\rho^{\mu_1}}
\int_{ (0,T)\cap (t_0-\rho;t_0+\rho)}
\left[
\sup_{%\substack
{x \in \R^n \atop  \rho\,>\,0 }}
\frac{1}{\rho^\lambda}
\int_{B_\rho(x) }
\left(  M_\eta(f)  \right)^{\frac{p}{\varepsilon}}(y,t)\,dy
\right]^{\frac{q_1}{q}}
dt
\leq \,
\end{equation}

\begin{equation}\label{63}
\leq \,
C\,
\frac{1}{\rho^{\mu_1}}
\int_{ (0,T)\cap (t_0-\rho;t_0+\rho)}
\left[
\sup_{%\substack
{x \in \R^n \atop  \rho\,>\,0 }}
\frac{1}{\rho^\lambda}
\int_{B_\rho(x) }
  |f|^p (y,t) \,dy
\right]^{\frac{q_1}{p}} dt
\end{equation}
the last term is less or equal than
\begin{equation}\label{64}
C\,
\sup_{%\substack
{t_0%, t
 \in (0,T)\atop  \rho>0 }}\,
\frac{1}{\rho^{\mu_1}}
\int_{ (0,T)\cap (t_0-\rho;t_0+\rho)}
\left[
\sup_{%\substack
{x \in \R^n \atop  \rho\,>\,0 }}
\frac{1}{\rho^\lambda}
\int_{B_\rho(x) }
  |f|^p (y,t) \,dy
\right]^{\frac{q_1}{p}} dt.
\end{equation}
Finally, we have
\begin{equation}\label{65}
\sup_{%\substack
{t_0 %, t
 \in (0,T)\atop  \rho>0 }}\,
\frac{1}{\rho^{\mu_1}}
\int_{ (0,T)\cap (t_0-\rho;t_0+\rho)}
\left[
\sup_{%\substack
{x \in \R^n \atop  \rho\,>\,0 }}
\frac{1}{\rho^\lambda}
\int_{B_\rho(x) }
  (M_\eta f)^q  (y,t)\,dy
\right]^{\frac{q_1}{q%p
}}
dt \leq
\end{equation}
\begin{equation}\label{66}
\leq
C\,
\sup_{%\substack
{t_0%, t
 \in (0,T)\atop  \rho>0 }}\,
\frac{1}{\rho^{\mu_1}}
\int_{ (0,T)\cap (t_0-\rho;t_0+\rho)}
\left[
\sup_{%\substack
{x \in \R^n \atop  \rho\,>\,0 }}
\frac{1}{\rho^\lambda}
\int_{B_\rho(x) }
  |f|^p (y,t) \,dy
\right]^{\frac{q_1}{p}}dt.
\end{equation}
Elevating both sides to $\frac{1}{q_1}, $ we have
\begin{equation}\label{67}
\| M_\eta\,f\|_{L^{q_1,\mu_1}(0,T,L^{q,\lambda}(\R^n))}
\leq \,C\,\,
\|f\|_{L^{q_1,\mu_1}(0,T,L^{p,\lambda}(\R^n))}.
\end{equation}

\subsection{%Applications to %Partial Differential Equations
Estimates of singular integral operators and commutators } %Section 5
Let $k(x,y) $ be a variable Calder\'on-Zygmund %variable
kernel for a.e. $x \in \R^{n+1}, $ %of mixed homogeneity,
 $f \in L^{q,\mu}(0,T,L^{p,\lambda}(\R^n)) $ with $ 1< p,q<\infty $ $0<\lambda, \mu <n, %\alpha
$   $a \in BMO(\R^{n+1}).$ For $\varepsilon >0$ let us define the operator $K_\varepsilon$ and the commutator $C_\varepsilon[a,f], $ as follows
\begin{equation}\label{67,1}
K_\varepsilon f(x) = \int_{\rho(x-y)>\varepsilon} %K
k(x,x-y) f(y) dy
\end{equation}

\begin{equation}\label{67,2}
C_\varepsilon[a,f]= K_\varepsilon (a f)(x) - a(x) K_\varepsilon f(x) =
\int_{\rho(x-y)>\varepsilon} %K
k(x,x-y) [a(x)- a(y)] f(y) dy.
\end{equation}
In the next theorem we prove that $K_\varepsilon f$ and $C_\varepsilon[a,f] $ are, uniformly  in $\varepsilon, $ bounded from $L^{q,\mu}(0,T,L^{p,\lambda}(\R^n))$ into itself. This fact allows us to let $\varepsilon \to 0 $ obtaining as limits in $L^{q,\mu}(0,T,L^{p,\lambda}(\R^n)) $ the following singular integral and commutator
\begin{equation}\label{67,3}
K f(x)= P.V. \int_{\R^n} k(x,x-y) f(y) dy = \lim_{\varepsilon \to 0} K_\varepsilon f(x)
\end{equation}
\begin{equation}\label{67,4}
C[a, f](x)= P.V. \int_{\R^n} k(x,x-y)[a(x)-a(y)] f(y) dy = \lim_{\varepsilon \to 0} C_\varepsilon [a,f](x)
\end{equation}

These operators are bounded in the class $L^{q,\mu}(0,T,L^{p,\lambda}(\R^n)). $
%fino a qui: [PalSoft 2004 Pot Analysis, righe 11-23]

\begin{thm}\label{thm5-3}\hspace*{-0.6em}\textbf{.}%(\cite{PS} Theorem 2.1 and Corollary 2.8$\,$)
%Pala-Soft Pot.An.2004, pg.244] e [DPR(1999)Theorem2.4]
Let $k(x,y) $ be a variable
Calder\'on-Zygmund kernel, for a.e. $x \in \R^{n+1}$, %of mixed homogeneity,
 $1<p,q<\infty, $  $0<\lambda,\mu <n%\alpha
  $ %($\alpha$ non n?)
   and $a \in VMO(\R^{n+1}).\,$

For any $f \in L^{q,\mu}(0,T,L^{p,\lambda}(\R^n))\,\,$ the singular integrals $K\,f,\,$ $C[a,f]\,\in L^{q,\mu}(0,T,L^{p,\lambda}(\R^n)).\,\,$ exist as limits in $L^{q,\mu}(0,T,L^{p,\lambda}(\R^n)), $ for $\varepsilon \to 0, $  of $K_\varepsilon f$ and $C_\varepsilon [a,f], $ respectively.
Then,  the operators $K\,f,\,C[a,f]\,: L^{q,\mu}(0,T,L^{p,\lambda}(\R^n))\,\to L^{q,\mu}(0,T,L^{p,\lambda}(\R^n))\,$ are bounded and satisfy the following inequalities

\begin{equation}\label{68}
\|Kf\|_{L^{q,\mu}(0,T,L^{p,\lambda}(\R^n))\,\,}\leq c \|f\|_{L^{q,\mu}(0,T,L^{p,\lambda}(\R^n))\,\,}
\end{equation}

\begin{equation}\label{69}
\|C[a,f]\|_{L^{q,\mu}(0,T,L^{p,\lambda}(\R^n))\,}\leq c \|a\|_* \|f\|_{L^{q,\mu}(0,T,L^{p,\lambda}(\R^n))\,}
\end{equation}

where $c=c(n,p,\lambda,\alpha,K), $ the dependence on $K $ is through the constant $c(\beta)$ in Definition \ref{Def5.2} part 2), for suitable $\beta. $

Moreover, for every $\epsilon >0 $ there exists $\rho_0>0 $ such that, if $B_r $ is a ball with radius $r\,$ such that $0<r<\rho_0,\,$ $k(x,y) $ satisfies the above assumptions %in $B_r $
 and $f \in L^{q,\mu}(0,T,L^{p,\lambda}(B_r)),$ we have
\begin{equation}\label{70}
\|C[a,f]\|_{L^{q,\mu}(0,T,L^{p,\lambda}(B_r))\,}\leq c \,\epsilon \,\|f\|_{L^{q,\mu}(0,T,L^{p,\lambda}(B_r))\,}
\end{equation}
for some constant %$c=c(n,p,\lambda,\M)$ being $M=\max_{|j|\leq 2n} \left|\frac{\partial^j}{\partial z^j} %k(x,y) \right|_{L^\infty(B \times \Sigma)}.$
$c$ independent of $f.$

\end{thm}

Proof. %We know from \cite{PS} that % Themtem 2.1 [Pal Soft 2004 Pot Analysis] o [Chiarenza-Frasca Rendiconti di mat., Theorem 3]
%\begin{equation}\label{71}
%\|Kf\|_{L^{p,\lambda}(\R^n)\,\,}\leq c \|f\|_{L^{p,\lambda}(\R^n)\,\,}
%\end{equation}
%then
%Let
For every $t \in (0,T),\,$  from the known inequality (see e.g. \cite{CF})

%(CORRETTO AGGIUNGENDO LA VARIABILE "$t$" ?)
\begin{equation}\label{72}
\sup_{%\substack
{x \in \R^n \atop  \rho\,>\,0 }}
\frac{1}{\rho^\lambda}
\int_{ B_\rho(x) }  |(K\,f)(y,t)|^p dy
\leq c
\sup_{%\substack
{x \in \R^n \atop  \rho\,>\,0 }}
\frac{1}{\rho^\lambda}
\int_{ B_\rho(x) } |f(y,t)|^p  \,dy,\quad  %\forall f \in L^p_\omega (\R^n)
\end{equation}
%true forall $t \in (t_0-\rho, t_0+\rho), $
 fixing $t_0 \in (0,T), \,$
 elevating to $\frac{q}{p}, $  integrating in $(0,T) \cap (t_0-\rho, t_0+\rho), $ multiplying for $\rho^{-\mu}, $% and taking the quii ,
 we have
\begin{equation}\label{73}
\frac{1}{\rho^\mu} \int_{(0,T) \cap (t_0-\rho, t_0+\rho)}
\left(
\sup_{%\substack
{x \in \R^n \atop  \rho\,>\,0 }}
\frac{1}{\rho^\lambda}
\int_{ B_\rho(x) }  |(K\,f)(y,t)|^p  dy
\right)^{\frac{q}{p}}
dt
\leq\,
\end{equation}
\begin{equation}\label{74}
\leq c
\frac{1}{\rho^\mu} \int_{(0,T) \cap (t_0-\rho, t_0+\rho)}
\left(
\sup_{%\substack
{x \in \R^n \atop  \rho\,>\,0 }}
\frac{1}{\rho^\lambda}
\int_{ B_\rho(x) } |f(y,t)|^p  \,dy
\right)^{\frac{q}{p}}
dt
\end{equation}
then, we have
\begin{equation}\label{75}
\sup_{%\substack
{t_0%, t
 \in (0,T)\atop  \rho>0 }}\,
\frac{1}{\rho^\mu} \int_{(0,T) \cap (t_0-\rho, t_0+\rho)}
\left(
\sup_{%\substack
{x \in \R^n \atop  \rho\,>\,0 }}
\frac{1}{\rho^\lambda}
\int_{ B_\rho(x) }  |(K\,f)(y,t)|^p  dy
\right)^{\frac{q}{p}}
dt
\leq\,
\end{equation}
\begin{equation}\label{76}
\leq\, c\, \sup_{%\substack
{t_0 %, t
 \in (0,T)\atop  \rho>0 }}\,
\frac{1}{\rho^\mu} \int_{(0,T) \cap (t_0-\rho, t_0+\rho)}
\left(
\sup_{%\substack
{x \in \R^n \atop  \rho\,>\,0 }}
\frac{1}{\rho^\lambda}
\int_{ B_\rho(x) }  |f(y,t)|^p  dy
\right)^{\frac{q}{p}}
dt
%\leq\,
\end{equation}
elevating to $\frac{1}{q}, $ we get the conclusion for $K\,f. $
Similar is the proof of (\ref{69}), starting from the inequality %(see \cite{PS}). %[Pal-Soft 2004 Pot Analysis, %Thm 2.1]
%Then, from \cite{DR1}  we have
\begin{equation}\label{77}
\|C[a,f]\|_{L^{p,\lambda}(\R^n)\,\,}\leq c \|a\|_*\|f\|_{L^{p,\lambda}(\R^n)\,\,}.
\end{equation}
Finally, using the $VMO$ assumption, if we %set that $\|a\|_*\,<\,<\,1\,$
fix $\rho_0$ such that $\eta(\rho_0)< \epsilon, $
we get the conclusion.
Let us remark that the result is also true if we assume $a$ defined only in some ball with $\|a||_*<\epsilon. $

\section{%Regularity Results
Applications to Partial Differential Equations } %Section 6
\setcounter{equation}{0}
\setcounter{thm}{0}
%See Pal Soft Potential Analysis 2004, pg.250

As application of the previous results we obtain a regularity result for strong solutions to the nondivergence form parabolic equations.

Precisely, let $n \geq 3, $  $Q_T%:
= \Omega' \times (0,T)\,$ be a cylinder of $\R^{n+%-
1
} \,$  of base $\Omega'\subset \R^{n%+1%-1
}. $
In the sequel let us set $ x=(x',t)=(x'_1,x_2,\ldots,x'_{n%-1
},t)\,$ a generic point in $Q_T, $
%$f \in L^p(Q_T)$ $ p(1;+\infty) $
$f \in L^{q,\mu}(0,T,L^{p,\lambda}(\Omega')), $ $1<p,q<\infty, $ $0<\lambda, \mu< n $ and
\begin{equation}\label{77,5}
L\,u=\,u_t\,-\, \Sigma_{i,j=1}^{n%-1
}\, a_{ij}(x',t)  \frac{\partial^2\,u}{\partial x'_i\partial x'_j}%=f(x',t) ,
\end{equation}
 where
\begin{equation}\label{78}
\,a_{ij}(x',t) = a_{ji}(x',t) , \qquad\forall i,j=1,\ldots,n, \qquad a.\,e.\, x \in Q_T \,%(Q_T ok ?) %x' \in \Omega' \,(\Omega\,or\, B ?)
\end{equation}
\begin{equation}\label{79}
\,\exists \nu > 0\,:\, \nu^{-1} |\xi|^2 \leq
\Sigma_{i,j=1}^n\, a_{ij}(x',t)\,\xi_i \xi_j \leq \nu |\xi|^2, \qquad a.\,e.\,\,\, in \,\,\,Q_T, \forall \xi \in \R^{n%-1
}
\end{equation}
\begin{equation}\label{80,1}
\,a_{ij}(%\cdot
x',t) \in
VMO( Q_T%\R^{n}
)\cap L^\infty(Q_T), \qquad \forall i,j=1,\ldots,n,.
\end{equation}
%La def. di strong solution e' presa da: [Pal Soft 2004 Potential Analysis, pg260]
Let us consider
\begin{equation}\label{80,2}
Lu(x',t)\,%=\,u_t\,-\, \Sigma_{i,j=1}^{n-1}\, a_{ij}(x',t)  \frac{\partial^2}{\partial x'_i\partial x'_j}
=f(x',t).
\end{equation}
A strong solution to (\ref{80,2}) is a function $u(x) \in L^{q,\mu}(0,T,L^{p,\lambda}(\Omega')) $ with all its weak derivatives $D_{x'_i}u, $ $D_{x'_i x'_j}u, $ $i, j=1,\ldots, n%-1
$ and $D_t u, $ satisfying (\ref{80,2})
, $\forall x \in Q_T. $

Let us now fix the coefficient $x_0=(x'_0, t_0) \in Q_T. $ and consider the fundamental solution of $L_0=L(x_0), $ is given, for $\tau>0, $ by
\begin{equation}\label{81}
\Gamma(x_0;\theta) = \Gamma (x'_0,t_0;\zeta,\tau)= \frac{(4\pi \tau)^{\frac{1-n}{2}}}{\sqrt{a^{ij}(x_0)} } exp\left(-
\frac{A^{ij}(x_0)\zeta_i\zeta_j }{4\,\tau}
\right)
\end{equation}
 that is equals to zero if $\tau \leq 0, $ being $A^{ij}(x_0)$ the entries of the inverse matrix $\{a^{ij}(x_0) \}^{-1}.$

The second order derivatives with respect to $\zeta_i$ and $\zeta_j,  $  denoted by  $\Gamma_{ij}(x_0,t_0;\zeta,\tau), $ $i,j=1,\ldots,n%-1
, $  and $\Gamma_{ij} (x;\theta), $ are kernels of mixed homogeneity in the sense that $\alpha_1=\ldots, \alpha_{n-1}=1, $ $\alpha_n=2$ (it follows that $\alpha= n+1$)

\begin{thm}\label{thm6-1}\hspace*{-0.6em}\textbf{.}%
Let $n \geq 3, a_{ij}\in VMO(Q_T)\cap L^\infty(Q_T),\,$ $B_r \subset\subset \Omega'$ a ball in $\R^{n%-1
}$ %(it is correct "a ball ?")

Then, for every $u$ having compact support in $B_r\times (0,T),\,$
solution of $Lu = f $ such that $D_{x'_i x'_j}u \in  L^{q,\mu}(0,T,L^{p,\lambda}(B_r)) \,$ $\forall i,j=1,\ldots, n%-1
, $ there exists $r_0=r_0(n.p,\nu, \eta%_\rho
)$ %(e' corretto $\eta_{\rho_0}$?)
 such that, if $r<r_0,$ then
\begin{equation}\label{82,1}
\|D_{x'_ix'_j}\,u \|_{L^{q,\mu}(0,T,L^{p,\lambda}(B_r))} \leq
C \|Lu\|_{L^{q,\mu}(0,T,L^{p,\lambda}(B_r))},\qquad i,j=1,\ldots,n%-1
\end{equation}
\begin{equation}\label{82,2}
\|\,u_t \|_{L^{q,\mu}(0,T,L^{p,\lambda}(B_r))} \leq
C \|Lu\|_{L^{q,\mu}(0,T,L^{p,\lambda}(B_r))},
\end{equation}
%being $u$ with compact support in $B_r \times (0,T) $
% or $u(\cdot, t) $ having compact support in $B_r.$
\end{thm}

Proof.
%Vd [Bramanti  Cerutti pg.1739]
Let $C_t =\{  v%u
\in C^\infty_0({\cal A}): v%u
(x',0)=0, {\cal A} =\R^{n+1} \cap\{t \geq 0\}\} $ and $u \in C_t. $
The local representation formula for the second order spatial derivatives
%, for any function
of $u $ %v \in C_t,%^\infty_0(B_r)$
  (see \cite{BC}),  is the following
 %(ha senso in $B_r$?)   $v(x',0)=0 $ is (see [Bramanti-Cerutti])
$$%\begin{equation}\label{83}
\!\!\!\!\!\!\!\!\!\!D_{x'_i x'_j}u(x)=%P.V.
 \lim_{\varepsilon \to 0 }\int_{\rho(x-y)>\varepsilon
 %B_r
 }\!\!\!\!\!\!\!\! \Gamma_{ij}(x;x-y) Lu(y) dy+%P.V.
%+\lim_{\varepsilon \to 0}  \int_{ \rho(x-y)>\varepsilon} \!\!\!\!\!\!\!\!\Gamma_{ij} (x;x-y) \Sigma_{h,k=1}^n  [a^{h %k}(y) - a^{h k}(x)] \cdot D_{y'_h\,y'_k} u(y) dy +
$$%\end{equation}
\begin{equation}\label{84}
+\lim_{\varepsilon \to 0}  \int_{ \rho(x-y)>\varepsilon} \!\!\!\!\!\!\!\!\Gamma_{ij} (x;x-y) \Sigma_{h,k=1}^n  [a^{h k}(y) - a^{h k}(x)] \cdot D_{y'_h\,y'_k} u(y) dy
+\,L(x) \int_\Sigma%{S^{n+1}}
 \nu_i(y) \Gamma_j(x;y) d\sigma,
\end{equation}
forall $i, j =1, \dots, n, $ and for $x$ in the support of $u, $ being %with
 $\nu_i(y) $ the i-th component of the unit outward normal to $\Sigma%S^{n+1}
  $ at $y \in %S^{n+1}
  \Sigma. $

From %Theorem \ref{thm5-3} we get the conclusion.
(\ref{68}) and (\ref{69}) we get the first inequality. Let us now observe that
\begin{equation}\label{85}
u_t= Lu + \Sigma_{i,j=1}^{n%-1
}\, a_{ij}(x',t)  \frac{\partial^2\,u}{\partial x'_i\partial x'_j}
\end{equation}
and the second inequality $(\ref{82,2}) $ is proved.

%{\bf Compliance with Ethical Standards}\\
%The author declare that she has no conflict of interest.

\footnotesize \noindent
{\sc Acknowledgments.} \
This work was conceived when the first author stayed at Mittag-Leffler Institute, in occasion of the research program  "Homogenization and Random Phenomenon". %, fall 2014.
 She wish to thank this institution for his hospitality.

% BIBLIOGRAFIA ---------------------------------------------------------
%\newpage

%{\bf References}
\end{document}